\documentclass{ws-procs9x6}
\usepackage{amsmath,amssymb,amscd,amsfonts,verbatim,graphicx}
\usepackage{xypic}

\newtheorem{lem}[theorem]{Lemma}
\newtheorem{propos}[theorem]{Proposition}

\newtheorem{cor}[theorem]{Corollary}

\newcommand{\co}{\colon\thinspace}    

\begin{document}

\title{Link groups and the A-B slice problem}

\author{Vyacheslav S. Krushkal}

\address{Department of Mathematics, University of Virginia, \\
Charlottesville, VA 22904-4137\\
Current address: Kavli Institute for Theoretical Physics\\
Santa Barbara, CA 93109\\
E-mail: krushkal@virginia.edu}

\centerline{\em Dedicated to the memory of Xiao-Song Lin}

\bigskip

\begin{abstract}
The $A-B$ slice problem is a reformulation of the topological
$4-$dimensional surgery conjecture in terms of
decompositions of the $4-$ball and link homotopy. We show that
{\em link groups}, a recently developed invariant of $4-$manifolds,
provide an obstruction for the class of {\em model decompositions},
introduced by M. Freedman and X.-S. Lin. This unifies and extends the
previously known partial obstructions in the $A-B$ slice program.
As a consequence, link groups satisfy Alexander duality when restricted
to the class of model decompositions, but not for general submanifolds
of the $4-$ball.
\end{abstract}


\bodymatter

\section{Introduction} \label{introduction}

The surgery conjecture, a core ingredient in the geometric
classification theory of topological $4-$manifolds, remains an open
problem for a large class of fundamental groups. The results to date
in the subject: the disk embedding conjecture, and its corollaries
-- surgery and s-cobordism theorems for {\em good} groups
\cite{Freedman0, FQ, FT1, KQ} -- show
similarities of classification of
topological $4-$manifolds with the theory in
higher dimensions. On the other hand, it has been
conjectured \cite{Freedman1} that surgery fails for (non-abelian) free fundamental
groups.

The $A-B$ slice problem \cite{F2} is a reformulation of the surgery
conjecture for free groups which seems most promising in terms of
the search for an obstruction. In this approach one considers
smooth codimension zero decompositions $D^4=A_i\cup B_i$ of the $4-$ball, extending
the standard genus one Heegaard decomposition of the $3$-sphere. (A precise definition
is given in section \ref{surgery overview}, also see figure \ref{fig:2D}.) Then the
problem is formulated in terms of the existence of disjoint embeddings of the submanifolds $A_i, B_i$
in $D^4$ with a prescribed homotopically essential link in $S^3=\partial D^4$
as the boundary condition. The central case corresponds to the link equal to the
Borromean rings.
The problem may be phrased in terms of the existence of a suitably formulated
non-abelian Alexander duality in dimension $4$. Recently this approach has been
sharpened and now there is a precise, axiomatic description of what
properties an obstruction, which in this context is an invariant of
decompositions of $D^4$, should satisfy.

The $A-B$ slice formulation of surgery was introduced by Freedman \cite{F2}
and further extensively studied by Freedman-Lin \cite{FL}. In particular, the latter paper
introduced a family of {\em model decompositions} which appear to approximate,
in a certain algebraic sense, an arbitrary decomposition $D^4=A\cup B$. This family of
decompositions is defined in section \ref{model}.
In this paper we use {\em link groups} of
$4-$manifolds, recently introduced by the author \cite{K},
to formulate an obstruction for the family of model decompositions:

\begin{theorem} \label{main theorem} \sl
Let $L$ be the Borromean rings, or more generally any homotopically essential link
in $S^3$. Then $L$ is not $A-B$ slice where each decomposition $D^4=A_i\cup B_i$
is a model decomposition.
\end{theorem}

The invariant using link groups formulated in the proof unifies and generalizes the
previously known partial obstructions\cite{FL, K1}
in the $A-B$ slice program. The definitions of link groups
and the underlying geometric notion of {\em Bing cells} are given in section \ref{definitions section}.

To place this result in the
geometric context of link homotopy, it is convenient to introduce
the notion of a robust $4-$manifold. Recall that a link $L$ in $S^3$ is
{\em homotopically trivial} \cite{M}  if its components bound disjoint
maps of disks in $D^4$. $L$ is called homotopically essential otherwise. (The Borromean
rings is a homotopically essential link with trivial linking numbers.)
Let $(M, {\gamma})$ be a pair
($4-$manifold, embedded curve in $\partial M$). The pair
$(M,{\gamma})$ is {\em robust} if whenever several copies $(M_i,
{\gamma}_i)$ are properly disjointly embedded in $(D^4, S^3)$, the
link formed by the curves $\{ {\gamma}_i\}$  in $S^3$ is
homotopically trivial. The following statement is a consequence of
the proof of theorem \ref{main theorem}:

\begin{cor} \sl
Let $D^4=A\cup B$ be a model decomposition. Then precisely one of the two
parts $A$, $B$ is  robust.
\end{cor}

It is interesting to note that there exist decompositions where {\em neither}
of the two sides is robust \cite{K2}.
The following question relates this notion to
the $A-B$ slice problem: {\sl given a decomposition $D^4=A\cup B$, is one of
the given embeddings $A\hookrightarrow D^4$, $B\hookrightarrow
D^4$ necessarily robust?}
(The definition of a robust embedding $e\co
(M,{\gamma})\hookrightarrow (D^4,S^3)$ is analogous to the
definition of a robust pair above, with the additional requirement
that each of the embeddings $(M_i, {\gamma}_i)\subset (D^4,S^3)$ is
equivalent to $e$.)

In a certain sense, one is looking in the $A-B$ slice
problem for an invariant of $4-$manifolds which is more flexible
than homotopy (so it satisfies a suitable version of Alexander
duality), yet it should be more robust than homology -- this is made
precise using Milnor's theory of link homotopy. The subtlety of the
problem is precisely in the interplay of these two requirements.
Following this imprecise analogy, we show that link groups
provide a step in construction of such a theory.


\section{Surgery
and the $A-B$ slice problem.} \label{surgery overview}

The $4-$dimensional topological surgery exact sequence (cf [FQ], Chapter 11),
as well as the $5-$dimensional
topological s-cobordism theorem, are known to hold for a class of
{\em good} fundamental groups. In the simply-connected case, this
followed from Freedman's disk embedding theorem \cite{Freedman0}
allowing one to represent hyperbolic pairs in ${\pi}_2(M^4)$ by
embedded spheres. Currently the class of good groups is known to
include the groups of subexponential growth \cite{FT1, KQ}
and it is closed under extensions and direct limits. There is a
specific conjecture for the failure of surgery for free groups
\cite{Freedman1}:

\begin{conjecture} \label{conjecture}
{\sl There does not exist a topological $4-$manifold
$M$, homotopy equivalent to ${\vee}^3 S^1$ and with $\partial M$
homeomorphic to ${\mathcal S}^0(Wh(Bor))$, the zero-framed surgery
on the Whitehead double of the Borromean rings.}
\end{conjecture}

In fact, this is one of a collection of canonical surgery problems
with free fundamental groups, and solving them is equivalent to the
unrestricted surgery theorem. The $A-B$ slice problem, introduced in
ref. \refcite{F2}, is a reformulation of the surgery conjecture, and it may
be roughly summarized as follows. Assuming on the contrary that the
manifold $M$ in the conjecture above exists, consider the
compactification of the universal cover $\widetilde M$, which is
homeomorphic to the $4-$ball \cite{F2}. The group of covering
transformations (the free group on three generators) acts on $D^4$
with a prescribed action on the boundary, and roughly speaking the
$A-B$ slice problem is a program for finding an obstruction to the
existence of such actions. Recall the definition of an $A-B$ slice
link \cite{F2, FL}.

\begin{figure}[ht]
\centering
\includegraphics[width=3.8cm]{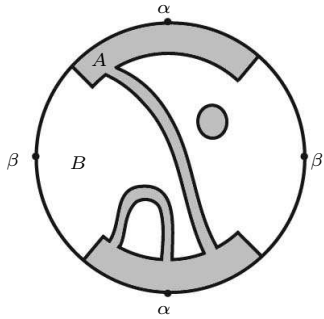}
{\scriptsize
    \put(-59,109){${\alpha}$}
    \put(-59,-5){${\alpha}$}
    \put(-84,89){$A$}
    \put(-116,51){${\beta}$}
    \put(-1,51){${\beta}$}
    \put(-92,50){$B$}}
    \vspace{.3cm}
\caption{A $2-$dimensional analogue of a decomposition $(A,{\alpha}),\, (B, {\beta})$:
$D^2=A\cup B$, $A$ is shaded; $({\alpha},\, {\beta})$ are linked $0-$spheres in $\partial D^2$.}
\label{fig:2D}
\end{figure}

\begin{definition} \label{decomposition}
A {\em decomposition} of $D^4$ is a pair of smooth compact
codimension $0$ submanifolds with boundary $A,B\subset D^4$,
satisfying conditions $(1)-(3)$ below. (Figure \ref{fig:2D} gives a $2-$dimensional example of
a decomposition.)
Denote $$\partial^{+}
A=\partial A\cap
\partial D^4, \; \; \partial^{+} B=\partial B\cap \partial D^4,\; \;
\partial A=\partial^{+} A\cup {\partial}^{-}A, \; \; \partial
B=\partial^{+} B\cup {\partial}^{-}B.$$
(1) $A\cup B=D^4$,\\
(2) $A\cap B=\partial^{-}A=\partial^{-}B,$ \\
(3) $S^3=\partial^{+}A\cup \partial^{+}B$ is the standard genus $1$
Heegaard decomposition of $S^3$.
\end{definition}
\begin{definition} \label{A-B slice}
Given an $n-$component link $L=(l_1,\ldots,l_n)\subset S^3$, let
$D(L)=(l_1,l'_1,\ldots, l_n,l'_n)$ denote the $2n-$component link
obtained by adding an untwisted parallel copy $L'$ to $L$. The link
$L$ is {\em $A-B$ slice} if there exist decompositions $(A_i, B_i),
i=1,\ldots, n$ of $D^4$ and self-homeomorphisms ${\alpha}_i,
{\beta}_i$ of $D^4$, $i=1,\ldots,n$ such that all sets in the
collection ${\alpha}_1 A_1, \ldots, {\alpha}_n A_n, {\beta}_1
B_1,\ldots, {\beta}_n B_n$ are disjoint and satisfy the boundary
data: ${\alpha}_i({\partial}^{+}A_i)$  is a tubular neighborhood of
$l_i$ and ${\beta}_i({\partial}^{+}B_i)$ is a tubular neighborhood
of $l'_i$, for each $i$.
\end{definition}

The surgery conjecture holds for all groups if and only if the
Borromean Rings (and the rest of the links in the canonical family
of links) are $A-B$ slice \cite{F2}. Conjecture \ref{conjecture}
above can therefore be reformulated as saying that the Borromean
Rings are not $A-B$ slice.

As an elementary example, note that if a link $L$ is $A-B$ slice
where for each $i$ the decomposition $D^4=A_i\cup B_i$ consists of
$A_i=2-$handle $D^2\times D^2$, and $B_i=$ the collar on
${\partial}^+ B_i$, then $L$ is actually {\em slice}.

Of course the Borromean Rings is not a slice (or homotopically
trivial) link. However to show that a link is not $A-B$ slice, one
needs to eliminate {\em all} choices for decompositions $(A_i,
B_i)$.

\section{Link groups and Bing cells.}
\label{definitions section}

In this section we recall the definition of Bing cells and link groups of
$4-$manifolds, denoted
${\lambda}(M^4)$, introduced in Ref. \refcite{K}, in order to formulate
the invariant $I_{\lambda}$ used in the proof of theorem \ref{main theorem}.
The definition is inductive.

\begin{definition} \label{model cell}
A {\em model Bing cell of height $1$} is a smooth $4$-manifold $C$
with boundary and with a specified attaching curve ${\gamma}\subset
\partial C$, defined as follows. Consider a planar surface $P$ with
$k+1$ boundary components ${\gamma}, {\alpha}_1,\ldots,{\alpha}_k$
($k\geq 0$), and set $\overline P=P\times D^2$. Let $L_1,\ldots,
L_k$ be a collection of links, $L_i\subset {\alpha}_i\times D^2$,
$i=1,\ldots,k$. Here for each $i$, $L_i$ is the (possibly iterated)
Bing double of the core ${\alpha}_i$. Then $C$ is obtained from
$\overline P$ by attaching zero-framed $2$-handles along the
components of $L_1\cup\ldots\cup L_k$.

The surface $S$ (and its thickening $\overline S$) will be referred
to at the {\em body} of $C$, and the $2$-handles are the {\em
handles} of $C$.

A {\em model Bing cell $C$ of height $h$} is obtained from a model
Bing cell of height $h-1$ by replacing its handles with Bing cells
of height one. The {\em body} of $C$ consists of all (thickenings
of) its surface stages, except for the handles.
\end{definition}

\begin{figure}[ht] \label{model figure1}
\centering
\includegraphics[width=6.25cm]{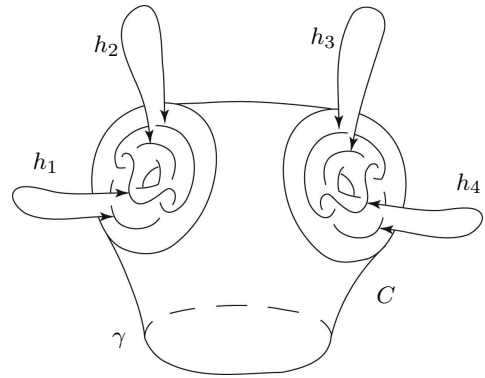}
{\small
    \put(-140,12){${\gamma}$}
    \put(-170,78){$h_1$}
    \put(-40,27){$C$}
    \put(-147,122){$h_2$}
    \put(-65,125){$h_3$}
    \put(-10,70){$h_4$}}
    \vspace{.3cm}
\caption{Example of a model Bing cell of height $1$: a schematic
picture}
\label{cell figure}
\end{figure}

\begin{figure}[ht] \label{model figure2}
\centering
\includegraphics[bb=0 -55 250 100, width=6.5cm]{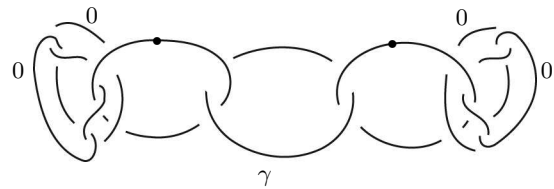}
{\small
    \put(-100,33){${\gamma}$}
    \put(-193,72){$0$}
    \put(-165,93){$0$}
    \put(-25,92){$0$}
    \put(7,72){$0$}}
    \put(-138,86){\circle*{3}}
    \put(-49,85){\circle*{3}}
    \vspace{.3cm}
\caption{A Kirby diagram of the model Bing cell in figure 1}
\label{cell figure Kirby}
\end{figure}

Figures \ref{cell figure}, \ref{cell figure Kirby} give
an example of a Bing cell of height $1$: a schematic
picture and a precise description in terms of a Kirby diagram. Here
$P$ is a pair of pants, and each link $L_i$ is the Bing double of
the core of the solid torus ${\alpha}_i\times D^2$, $i=1,2$.

\begin{remark}
To avoid a technical discussion, the definition presented here
involves only the links $L$ which are Bing doubles. To reflect this
difference, we reserve for these objects the term {\em Bing cells}
rather than the more general {\em flexible cells} discussed in Ref. \refcite{K}.
The definition in Ref. \refcite{K} involves more general
homotopically essential links, however just the Bing doubles suffice
for the applications in this paper.
\end{remark}

{\em Bing cells in a $4-$manifold} $M$ are defined as maps of model
Bing cells in $M$, subject to certain crucial disjointness
requirements. (In particular, this will be important for the
discussion of model decompositions in section \ref{model}.) Roughly
speaking, objects attached to different components of any given link
$L_i$ in the definition are required to be disjoint in $M$. To
formulate this condition rigorously, recall the definition of the
tree associated to a given Bing cell.

\subsection{The associated tree} \label{associated tree}
Given a Bing cell $C$, define the tree $T_C$ inductively: suppose
$C$ has height $1$. Then assign to the body surface $P$ (say with
$k+1$ boundary components) of $C$ the cone $T_P$ on $k+1$ points.
Consider the vertex corresponding to the attaching circle ${\gamma}$
of $C$ as the root of $T_P$, and the other $k$ vertices as the
leaves of $T_P$. For each handle of $C$ attach an edge to the
corresponding leaf of $T_P$. The leaves of the resulting tree $T_C$
are in $1-1$ correspondence with the handles of $C$.

\begin{figure}[t] \label{cell and tree}
\centering
\includegraphics[width=3.5cm]{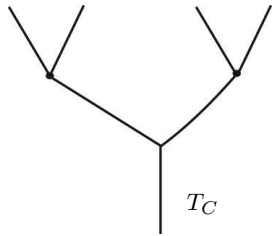}
\hspace{.7cm} {\small
    \put(-32,9){$T_C$}}
    \vspace{.3cm}
    \caption{The tree $T_C$ associated to the Bing cell $C$ in figures
    \ref{cell figure}, \ref{cell figure Kirby}.}
\label{tree figure}
\end{figure}

Suppose $C$ has height $h>1$, then it is obtained from a Bing cell
$C'$ of height $h-1$ by replacing the handles of $C'$ with Bing
cells $\{ C_i\}$ of height $1$. Assuming inductively that $T_{C'}$
is defined, one gets $T_C$ by replacing the edges of $T_{C'}$
associated to the handles of $C'$ with the trees corresponding to
$\{ C_i\}$. Figure \ref{tree figure} shows the tree associated to the Bing cell in
figure \ref{cell figure}.

Divide the vertices of $T_C$ into two types: the vertices (``cone points'')
corresponding to body (planar) surfaces are {\em unmarked}; the rest
of the vertices are {\em marked}. Therefore the valence of an
unmarked vertex equals the number of boundary components of the
corresponding planar surface. The marked vertices are in $1-1$
correspondence with the links $L$ defining $C$, and the valence of a
marked vertex is the number of components of $L$ plus $1$. It is
convenient to consider the $1-$valent vertices of $T_C$: its root
and leaves (corresponding to the handles of $C$) as unmarked. This
terminology is useful in defining the maps of Bing cells below. The
height of a Bing cell $C$ may be read off from $T_C$ as the maximal
number of marked vertices along a geodesic joining a leaf of $T_C$
to its root, where the maximum is taken over the leaves of $T_C$.

\begin{definition} \label{cell}
A {\em Bing cell} is a model Bing cell with a finite number of
self-plumbings and plumbings among the handles and body surfaces of
$C$, subject to the following disjointness requirement:

$\bullet$ Consider two surfaces $A,B$ (they could be handles or body
stages) of $C$. Let $a,b$ be the corresponding vertices in $T_C$.
(For body surfaces this is the corresponding unmarked cone point,
for handles this is the associated leaf.) Consider the geodesic
joining $a,b$ in $T_C$, and look at its vertex $c$ closest to the
root of $T_C$ -- in other words, $c$ is the first common ancestor of
$a,b$. If $c$ is a marked vertex then $A, B$ are required to be disjoint.

\medskip

In particular, self-plumbings of any handle and body surface are
allowed. In the example shown in figures 1, 2 above, the handle
$h_1$ is required to be disjoint from $h_2$, $h_3$ is disjoint from
$h_4$; all other intersections are allowed.

A {\em Bing cell in a $4$-manifold $M$} is an embedding of a Bing
cell into $M$. We say that its image is a realization of $C$ in $M$,
and abusing the notation we denote its image in $M$ also by $C$.
\end{definition}

The main technical result of Ref. \refcite{K} shows how Bing cells fit in
the context of Milnor's theory of link homotopy. This theorem is
used in the analysis of the invariant $I_{\lambda}$ below.

\begin{theorem} \label{homotopy}
{\sl If the components of a link $L\subset S^3=\partial D^4$ bound
disjoint Bing cells in $D^4$ then $L$ is homotopically trivial.}
\end{theorem}

Recall \cite{M} that a link $L$ in $S^3$ is {\em homotopically
trivial} if $L$ is homotopic to the unlink, so that different
components stay disjoint during the homotopy. The theorem above
builds on a classical result that if the components of $L$ bound
disjoint maps of {\em disks} in $D^4$ then $L$ is homotopically
trivial. The proof of theorem \ref{homotopy} is substantially more
involved than the argument in the classical case. This is due to the
topology of Bing cells which forces additional relations in the
fundamental group of the complement. The main new technical
ingredients in the proof are the {\em generalized Milnor group} and
an obstruction which is well-defined in the presence of this
additional indeterminacy \cite{K}.

The {\em link groups} ${\lambda}_n(M)$ are defined as $\{$based loops
in a $4-$manifold $M\}$ modulo loops bounding  Bing cells of height
$n$. These groups fit in a sequence of surjections
$${\pi}_1(M)\longrightarrow
{\lambda}_1(M)\longrightarrow{\lambda}_2(M)\longrightarrow\ldots$$

The groups ${\lambda}_n(M)$ are topological but
not in general homotopy invariants of $M$. In particular, they are
not correlated with the first homology $H_1(M)$, or more generally
with the quotients of ${\pi}_1(M)$ by the terms of its lower central
or derived series. Define ${\lambda}(M)$ to be the direct limit of
${\lambda}_n(M)$. Given a pair $(M,{\gamma})$ where $M$ is a $4-$manifold
and $\gamma$ is a specified curve in $\partial M$, consider the
invariant $I_{\lambda}(M,{\gamma})\in\{ 0,1\}$:
$$I_{\lambda}(M,{\gamma})=1\; \; {\rm if}\; \;
{\gamma}=1\in {\lambda}(M),$$ set
$I_{\lambda}(M,{\gamma})=0$ otherwise. When the choice of the attaching circle
$\gamma$ of $M$ is clear, we will abbreviate the notation to $I_{\lambda}(M)$.

\begin{remark} \label{gropes} For the interested reader we point out the
``geometric duality'' between Bing cells and gropes. Recall the
definition \cite{FQ}: A {\it grope} is a special pair (2-complex,
circle).  A grope has a {\em class} $k=1, 2,\dots, \infty$. For
$k=2$ a grope is a compact oriented surface $\Sigma$ with a single
boundary component. For $k>2$ a $k$-{\it grope} is defined
inductively as follow: Let $\{\alpha_i, \beta_i, i=1, \dots,
\rm{genus}\}$ be a standard symplectic basis of circles for
$\Sigma$.  For any positive integers $p_i, q_i$ with $p_i+q_i\ge k$
and $p_{i_0} + q_{i_0} = k$ for at least one index $i_0$, a
$k$-grope is formed by gluing $p_i$-gropes to each $\alpha_i$ and
$q_i$-gropes to each $\beta_i$. A grope has a standard,
``untwisted'' $4-$dimensional thickening, obtained by embedding it
into ${\mathbb R}^3$, times $I$.

Consider a more general collection of $2$-complexes, where at each
stage one is allowed to attach several parallel copies of surfaces.
Then one checks using Kirby calculus that model Bing cells are
precisely complements in $D^4$ of standard embeddings of such
generalized gropes. This observation is helpful in the analysis of
the $A-B$ slice problem, where gropes play an important role, see section
\ref{model}.
\end{remark}

\section{An obstruction for model decompositions.} \label{model}

In this section we show that the invariant $I_{\lambda}$ defined above
provides an obstruction for the family of model
decompositions. We start the proof of theorem \ref{main theorem} by
constructing the relevant decompositions of $D^4$. The simplest
decomposition $D^4=A\cup B$ where $A$ is the $2-$handle $D^2\times
D^2$ and $B$ is just the collar on its attaching curve, was
discussed in the introduction. Now consider the genus one surface
$S$ with a single boundary component ${\alpha}$, and set $A_1=S\times D^2$.
Moreover, one has to specify its embedding into
$D^4$ to determine the complementary side, $B$. Consider the
standard embedding (take an embedding of the surface in $S^3$, push
it into the $4-$ball and take a regular neighborhood.) Note that
given any decomposition, by Alexander duality the attaching curve of
exactly one of the two sides  vanishes in it homologically, at least
rationally. Therefore the decomposition $D^4=A_1\cup B_1$ may
be viewed as the first level of an ``algebraic approximation'' to an
arbitrary decomposition. The general {\em model decomposition of height $1$}
is analogous to the decomposition $D^4=A_1\cup B_1$, except that the surface
$S$ may have a higher genus.

\begin{figure}[ht]
\centering
\includegraphics[width=2.8cm]{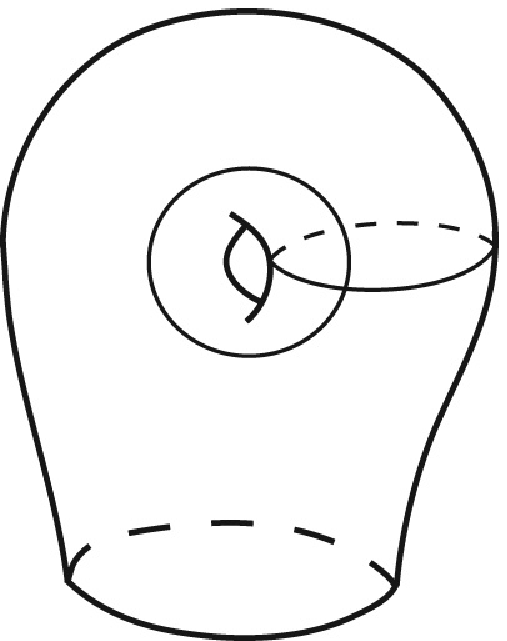} \hspace{3cm} \includegraphics[width=3.8cm]{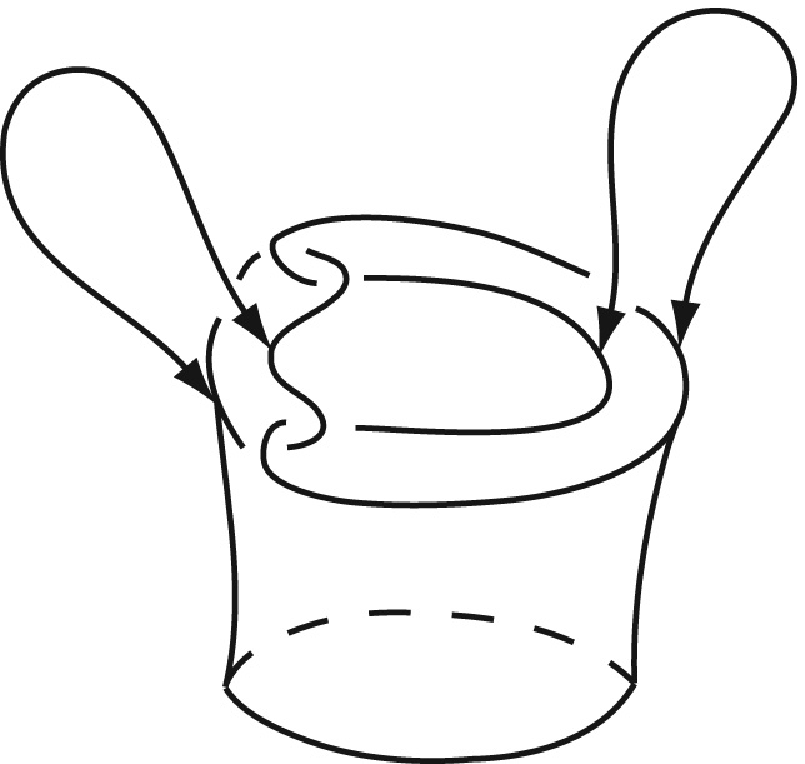}
{\small
    \put(-297,27){$A_1$}
    \put(-264,41){${\alpha}_1$}
    \put(-221,46){${\alpha}_2$}
    \put(-270,-7){${\alpha}$}
    \put(-80,-7){${\beta}$}
    \put(-126,73){$H_1$}
    \put(2,80){$H_2$}
    \put(-8,23){$B_1$}}
    \vspace{.3cm}
\caption{A model decomposition $D^4=A_1\cup B_1$ of height 1: a schematic (spine)
picture (figure 5) and a precise description in terms of Kirby diagrams, figure 6.}
\label{fig:surface}
\end{figure}

\begin{figure}[ht]
\centering
\includegraphics[width=4cm]{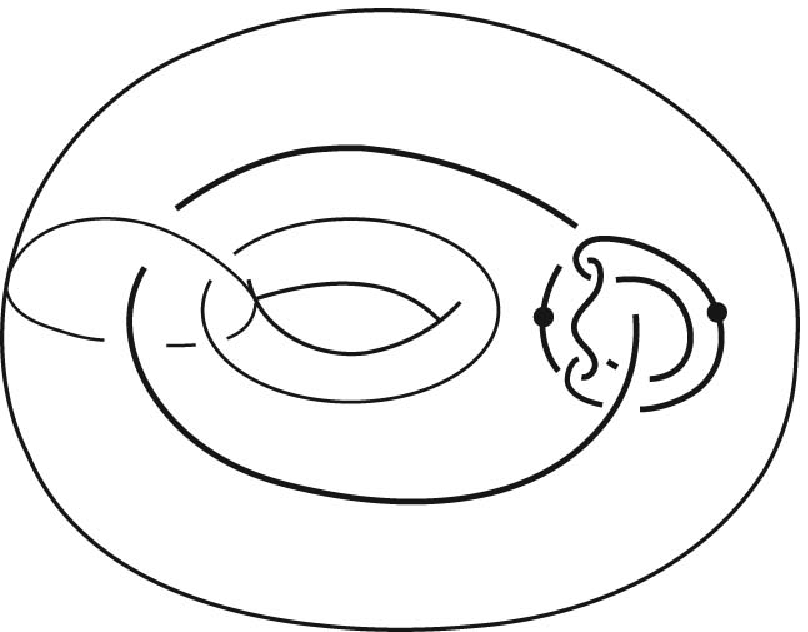} \hspace{2cm} \includegraphics[width=4cm]{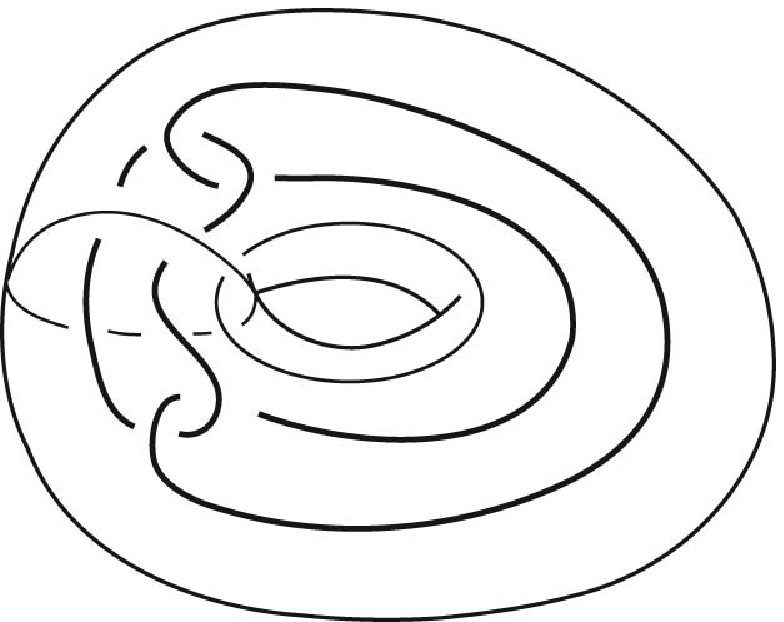}
{\small
    \put(-296,5){$A_1$}
    \put(-118,5){$B_1$}}
{\scriptsize
    \put(-270,20){$0$}
    \put(-230,29){$\alpha$}
    \put(-107,30){$0$}
    \put(-17,22){$0$}
    \put(-46,35){$\beta$}}
\caption{}
\label{fig:surface Kirby}
\end{figure}

\begin{propos} \label{surface complement}
\sl Let $A_1=S\times D^2$, where $S$ is the genus one
surface with a single boundary component $\alpha$. Consider the
standard embedding $(A_1, {\alpha}\times\{ 0\})\subset (D^4,
S^3)$. Then the complement $B_1$ is obtained from the collar
on its attaching curve, $S^1\times D^2\times I$, by attaching a pair
of zero-framed $2-$handles to the Bing double of the core of the
solid torus $S^1\times D^2\times\{1\}$, figures \ref{fig:surface}, \ref{fig:surface Kirby}.
\end{propos}

The proof is a standard exercise in Kirby calculus, see for example
Ref. \refcite{FL}. A precise description of these $4-$manifolds is given in
terms of Kirby diagrams in figure \ref{fig:surface Kirby}. Rather than considering
handle diagrams in the $3-$sphere, it is convenient to draw them in
the solid torus, so the $4-$manifolds are obtained from $S^1\times
D^2\times I$ by attaching the $1-$ and $2-$handles as shown in the
diagrams. To make sense of the ``zero framing'' of curves which are
not null-homologous in the solid torus, recall that the solid torus
is embedded into $S^3=\partial D^4$ as the attaching region of a
$4-$manifold, and the $2-$handle framings are defined using this
embedding.

This example illustrates the general principle that (in all examples
considered in this paper) the $1-$handles of each side are in
one-to-one correspondence with the $2-$handles of the complement.
This is true since the embeddings in $D^4$ considered here are all
standard, and in particular each $2-$handle is unknotted in $D^4$.
The statement follows from the fact that $1-$handles may be viewed
as standard $2-$handles removed from a collar, a standard technique
in Kirby calculus (see Chapter 1 in Ref. \refcite{Ki}.)
Moreover, in each of our examples the attaching
curve ${\alpha}$ on the $A-$side bounds a surface in $A$, so it has
a zero framed $2-$handle attached to the core of the solid torus. On
the $3-$manifold level, the zero surgery on this core transforms the
solid torus corresponding to $A$ into the solid torus corresponding
to $B$. The Kirby diagram for $B$ is obtained by taking the diagram
for $A$, performing the surgery as above, and replacing all zeroes
with dots, and conversely all dots with zeroes. (Note that the
$2-$handles in all our examples are zero-framed.)

Note that a distinguished pair of curves ${\alpha}_1, {\alpha}_2$,
forming a symplectic basis in the surface $S$, is determined as the
meridians (linking circles) to the cores of the $2-$handles $H_1,
H_2$ of $B_1$ in $D^4$. In other words, ${\alpha}_1$,
${\alpha}_2$ are fibers of the circle normal bundles over the cores
of $H_1, H_2$ in $D^4$.

\begin{figure}[ht]
\centering
\includegraphics[width=9cm]{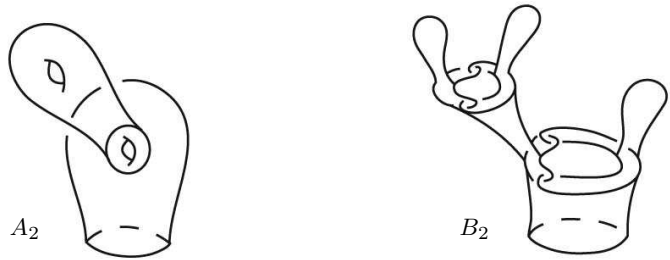}
\hspace{.7cm} {\small
    \put(-250,10){$A_2$}
    \put(-80,10){$B_2$}}
    \vspace{.3cm}
    \caption{A model decomposition $D^4=A_2\cup B_2$ of height $2$.}
\label{fig:models of height 2}
\end{figure}

An important observation \cite{FL} is that this construction may be
iterated: consider the $2-$handle $H_1$ in place of the original
$4-$ball. The pair of curves (${\alpha}_1$, the attaching circle
${\beta}_1$ of $H_1$) form the Hopf link in the boundary of $H_1$.
As discussed in the beginning of this section, it is natural to consider
two possibilities: either ${\alpha}_1$ or ${\beta}_1$ bounds a surface
in $H_1$. For simplicity of exposition, we again assume at this point
that this is a surface of genus one. The first possibility (${\alpha}_1$ bounds)
is shown in figure \ref{fig:models of height 2}: note that in this decomposition one
side, $A_2$, is a grope of
height $3$ (discussed in remark \ref{gropes}) and its complement $B_2$ is an
example of a Bing cell.

\begin{figure}[ht]
\centering
\includegraphics[width=3.6cm]{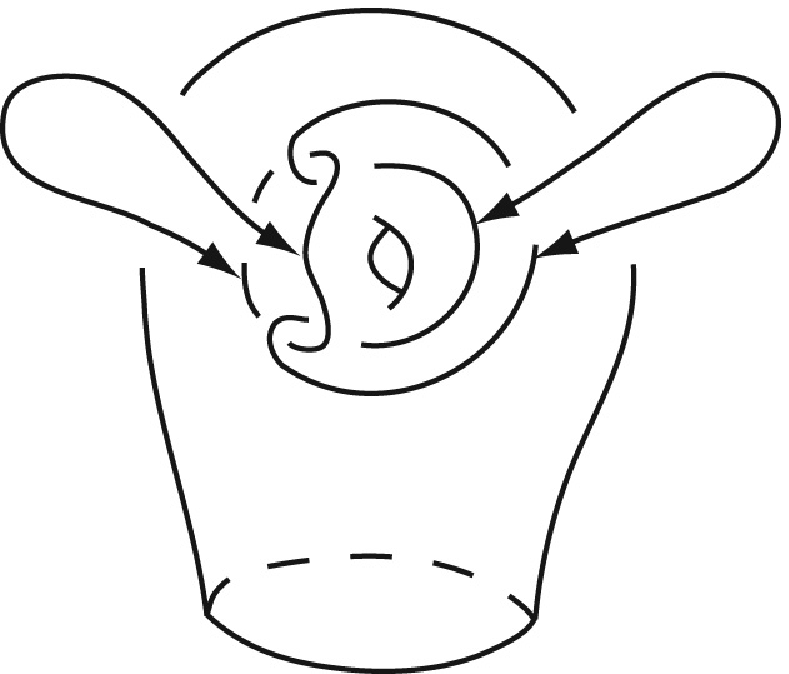}
 \hspace{1.5cm} \includegraphics[width=5.6cm]{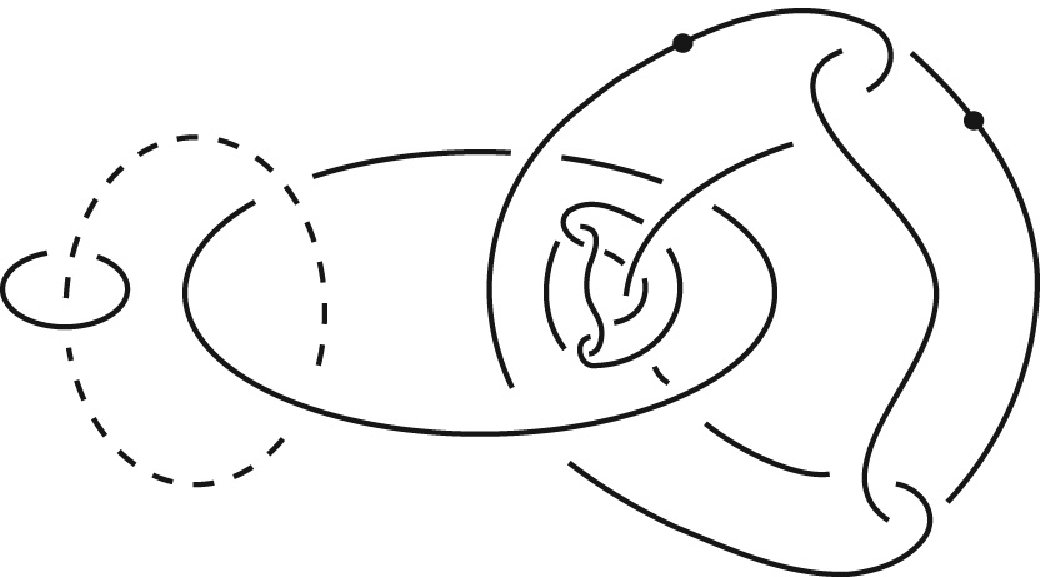}
{\small
    \put(-307,36){$A'_{2}$}
    \put(-290,-2){${\alpha}$}
    \put(-168,42){${\alpha}$}}
{\scriptsize
    \put(-99,13){$0$}
    \put(-81,42){$0$}
    \put(-53,42){$0$}}
\caption{}
\label{fig:A2}
\end{figure}

\begin{figure}[ht]
\centering
\includegraphics[width=3.5cm]{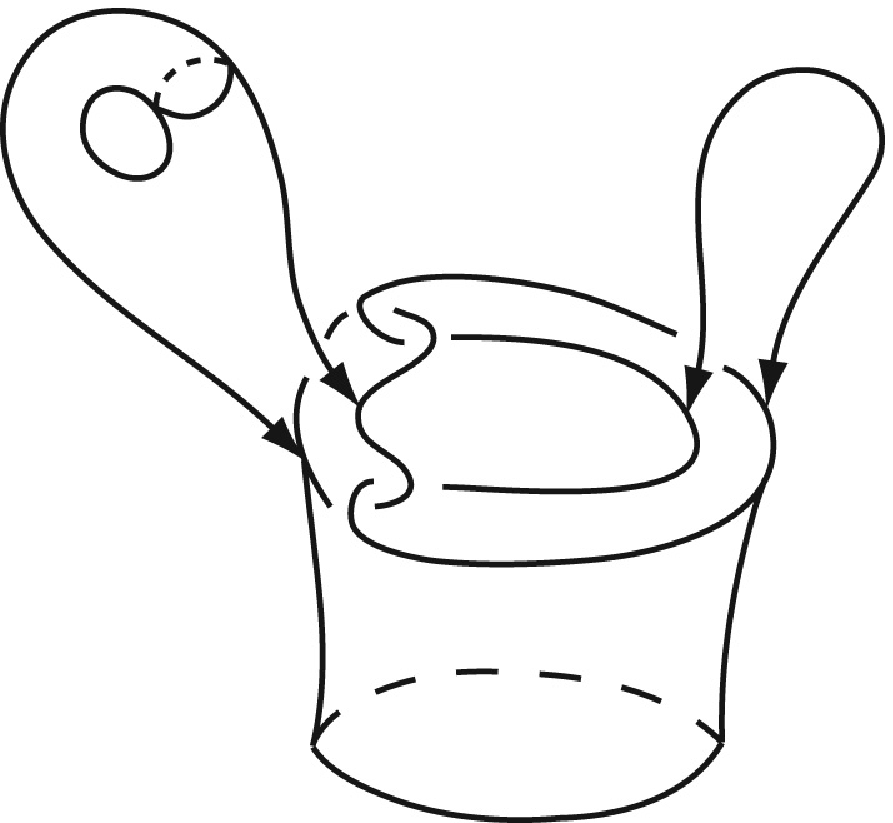} \hspace{1.8cm} \includegraphics[width=5.4cm]{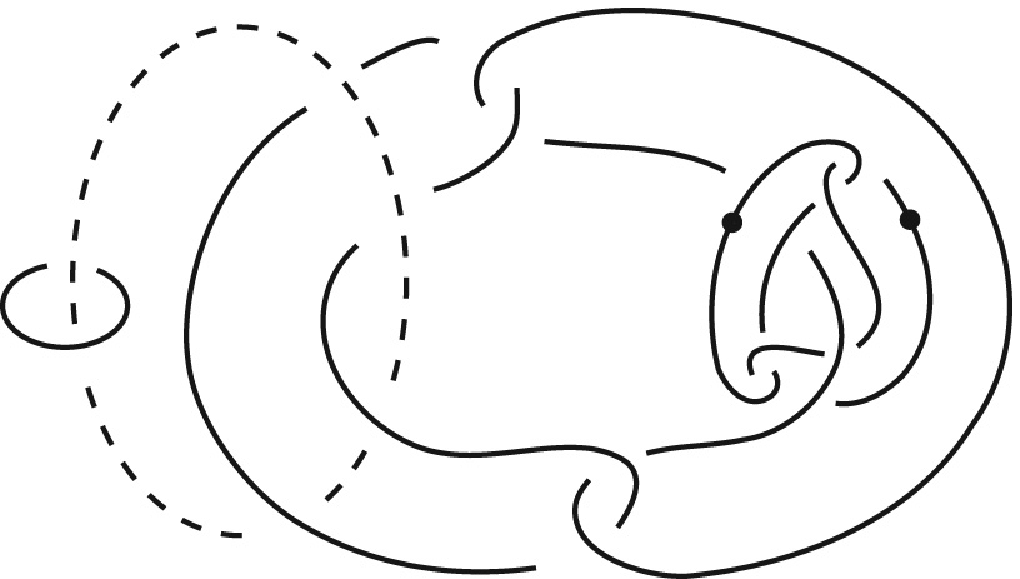}
{\small
    \put(-307,28){$B'_{2}$}
    \put(-290,0){${\beta}$}
    \put(-165,38){${\beta}$}}
{\scriptsize
    \put(-87,-9){$0$}
    \put(-47,-9){$0$}}
    \vspace{.3cm}
    \caption{Another example of a model decomposition
    $D^4=A'_2\cup B'_2$ of height $2$.}
\label{fig:A2 Kirby}
\end{figure}

\begin{figure}[ht]
\centering
\includegraphics[width=8cm]{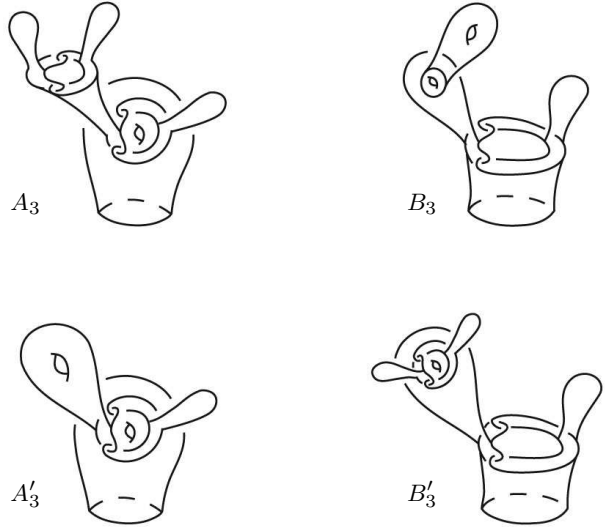}
{\small
    \put(-225,120){$A_3$}
    \put(-75,120){$B_3$}
    \put(-225,10){$A'_3$}
    \put(-75,10){$B'_3$}}
    \vspace{.3cm}
    \caption{Examples of model decompositions $D^4=A_3\cup B_3$, $D^4=A'_3\cup B'_3$
    of height $3$.}
\label{fig:models of height 3}
\end{figure}

Consider the second possibility: ${\beta}_1$ bounds a surface in $H_1$.
As discussed above, its complement in $H_1$ is given by two
zero-framed $2-$handles attached to the Bing double of ${\alpha}_1$.
Assembling this data, consider the new decomposition $D^4=A'_2\cup
B'_2$, figures \ref{fig:A2}, \ref{fig:A2 Kirby}. As above, the diagrams are drawn in solid tori
(complements in $S^3$ of unknotted circles drawn dashed in the
figures.) The decompositions $D^4=A_2\cup B_2$, $D^4=A'_2\cup B'_2$ are examples of model
decompositions of height $2$. To get a general decomposition of this type,
one also considers the alternative as above for the pair of curves ${\alpha}_2$,
${\beta}_2$ in the $4-$ball $H_2$. For simplicity of illustration, in the examples shown in
figures \ref{fig:models of height 2}
- \ref{fig:A2 Kirby} the curve ${\beta}_2$ bounds a surface of genus zero. One gets models of an
arbitrary height by an iterated application of the
construction above, and in general one considers (orientable) surfaces of an arbitrary genus at
each stage.  See figure \ref{fig:models of height 3}
for examples of model decompositions of height $3$.

It follows from theorem \ref{homotopy} that the following lemma implies our
main result, theorem \ref{main theorem}:

\begin{lem} \label{model lemma} \sl Let $D^4=A\cup B$ be a
model decomposition.
Then $$I_{\lambda}(A,{\alpha})+I_{\lambda}(B, {\beta})=1.$$
\end{lem}

Indeed, suppose a link $L=(l_1,\ldots, l_n)$ is $A-B$ slice where
each decomposition $D^4=A_i\cup B_i$, $i=1,\ldots,n$ is a model decomposition.
According to lemma \ref{model lemma}, the invariant $I_{\lambda}$ of precisely
one part of the decomposition equals $1$. For each $i$, denote $C_i=A_i$ if
$I_{\lambda}(A_i)=1$ and $C_i=B_i$ otherwise. Let ${\gamma}_i$ denote
the attaching curve of $C_i$. It follows from the definition  of $I_{\lambda}$
that ${\gamma}_i$ bounds a Bing cell in $C_i$. Since the collections $\{ {\alpha}_i\},
\{ {\beta}_i\}$ form the link $L$ and its parallel copy, the collection of curves
$({\gamma}_1,\ldots, {\gamma}_n$) is isotopic to $L$. This contradicts theorem \ref{homotopy}
since $L$ is homotopically essential. This concludes the proof of theorem \ref{main theorem}, assuming
lemma \ref{model lemma}.

\medskip

\noindent
{\em Proof of lemma \ref{model lemma}.} It suffices to prove that
given a model decomposition $D^4=A\cup B$, either ${\alpha}=1\in
{\lambda}(A)$ or ${\beta}=1\in {\lambda}(B)$. Then theorem
\ref{homotopy} implies that precisely one of these two possibilities
holds. The proof of the statement above is inductive. Given a model
decomposition of height $1$ (figure \ref{fig:surface}), observe that
one of the two parts of the decomposition - the handlebody $B_1$ in the
example in figure \ref{fig:surface} - is a model Bing cell of height $1$. (In this
case the planar surface $C$ in definition \ref{model cell} is the annulus.)
Therefore ${\beta}=1\in {\lambda}(B_1)$. In the case that $A_1$ is a surface
of genus $g>1$, the handlebody description of $B_1$ consists of first taking
$g$ parallel copies of the core curve of the solid torus, Bing doubling them and
then attaching zero-framed $2-$handles to the resulting link. One observes that
the attaching curve $\beta$ still bounds a model Bing cell of height $1$ in this
handlebody, indeed there are $g$ choices of Bing cells bounded by $\beta$.

\begin{figure}[ht]
\centering
\includegraphics[width=9cm]{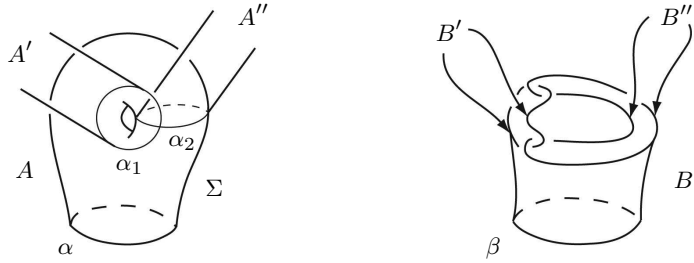}
{\small
    \put(-257,27){$A$}
    \put(-219,30){${\alpha}_1$}
    \put(-199,39){${\alpha}_2$}
    \put(-241,-2){${\alpha}$}
    \put(-260,71){$A'$}
    \put(-173,85){$A''$}
    \put(-185,20){$\Sigma$}
    \put(-79,-2){${\beta}$}
    \put(-98,79){$B'$}
    \put(-13,85){$B''$}
    \put(-8,23){$B$}}
    \vspace{.3cm}
\caption{Proof of lemma \ref{model lemma}: the inductive step.}
\label{fig:inductive}
\end{figure}

Suppose lemma is proved for model decompositions of height $\leq n$, and let
$D^4=A\cup B$ be a model decomposition of height $n+1$. The attaching curve
of either $A$ or $B$ is trivial in its first homology group. To be specific,
assume ${\alpha}=0\in H_1(A;{\mathbb Z})$. First assume the surface $\Sigma$
bounded by ${\alpha}$ has genus $1$. Then $A$ is obtained by attaching models
$A'$, $A''$ of height $\leq n$ to a symplectic basis of curves ${\alpha}_1, {\alpha}_2$
of $\Sigma$, figure \ref{fig:inductive}. Similarly, using the notation of figure \ref{fig:surface},
$B$ is obtained from the model $B_1$ of height $1$ by replacing its $2-$handles
$H_1, H_2$ by two models $B', B''$ of height $\leq n$. Here $D^4=A'\cup B', D^4=A''\cup B''$
are two decompositions for which lemma holds according to the inductive assumption.
Therefore $I_{\lambda}(A')+I_{\lambda}(B')=I_{\lambda}(A'')+I_{\lambda}(B'')=1$. Consider two
cases:

Case 1: $I_{\lambda}(B')=I_{\lambda}(B'')=1$

Case 2: At least one of $I_{\lambda}(A')$, $I_{\lambda}(A'')$ equals $1$.

\begin{figure}[ht]
\centering
\includegraphics[width=9cm]{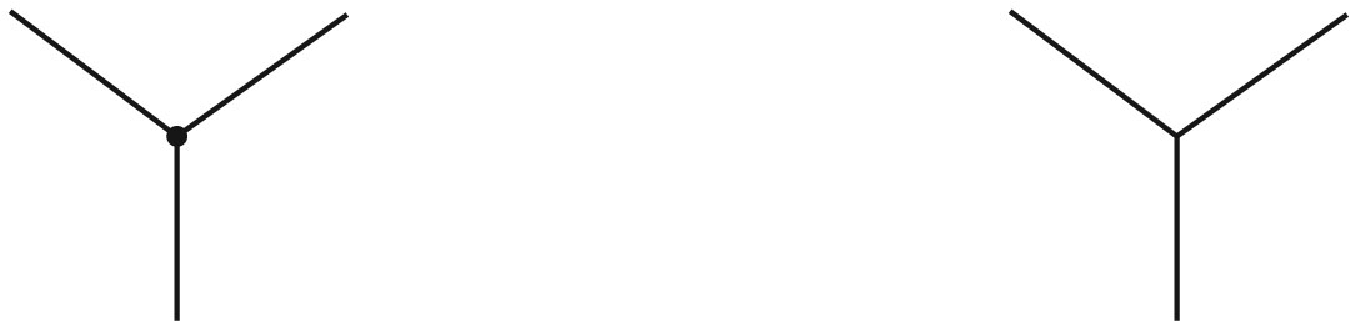}
{\small
    \put(-245,2){$T_C$}
    \put(-270,65){$T_{C'}$}
    \put(-186,65){$T_{C''}$}
    \put(-54,2){$T_C$}
    \put(-82,65){$T_{C'}$}
    \put(2,65){$T_{\overline C'}$}}
\caption{}
\label{fig:tree1}
\end{figure}

We claim that in the first case $I_{\lambda}(B)=1$ and in the second case
$I_{\lambda}(A)=1$. Consider case 1. By assumption, the attaching curve
${\beta}'$ of $B'$ bounds a Bing cell $C'$ in $B'$, and similarly the attaching curve
${\beta}''$ bounds a Bing cell $C''$ in $B''$. Consider the handlebody $C$ obtained from
$S^1\times D^2\times I$ by attaching $C', C''$ to the Bing double of the core of the solid
torus. The associated tree $T_C$ is illustrated on the left in figure \ref{fig:tree1}. (Note
that the trees $T_{C'}$, $T_{C''}$ join in a {\em marked} vertex.) Since
$B'$ and $B''$ are disjoint, there are no $C'-C''$ intersections. (Note that such intersections
are not allowed in the definition \ref{cell} of a Bing cell.)
Therefore the attaching curve $\beta$ bounds a Bing cell in $B$, and $I_{\lambda}(B,{\beta})=1$.

Consider the second case. Without loss of generality assume $I_{\lambda}(A')=1$, so
${\alpha}_1$ bounds a Bing cell $C'$ in $A'$. Surger the first stage surface $\Sigma$ along
${\alpha}_1$, the result is a pair of pants whose boundary consists of ${\alpha}$ and two copies
of ${\alpha}_1$. Consider two copies of $C'$ (denote them by $C'$ and $\overline C'$)
and perturb them so there are only finitely many intersections between surfaces in $C'$ and surfaces in $\overline C'$.
Consider the handlebody $C$ assembled from the (pair of pants)$\times D^2$ with $C', \overline C'$ attached
to it. The tree $T_C$ associated to $C$ is shown on the right in figure \ref{fig:tree1}; observe that
the trees $T_{C'}$, $T_{\overline C'}$ join in an {\em unmarked} vertex. Note that all
$C'-\overline C'$ intersections are of the type allowed in definition \ref{cell}, therefore
$\alpha$ bounds a Bing cell in $A$, and $I_{\lambda}(A,{\alpha})=1$.

In the case when the surface ${\Sigma}$ has genus $g>1$ the proof is analogous to the genus one case
discussed above. Specifically, $A$ is obtained by attaching models $A'_i$, $A''_i$, $i=1,\ldots, g$
to a symplectic basis of curves in $\Sigma$. The complements are denoted $B'_i, B''_i$. One observes that
if there exists $1\leq i\leq g$ such that
$I_{\lambda}(B'_i)=I_{\lambda}(B''_i)=1$, then $I_{\lambda}(B)=1$.
On the other hand, if for each $i$ either $I_{\lambda}(A'_i)$ or $I_{\lambda}(A''_i)$ equals $1$, then
$I_{\lambda}(A)=1$.
This concludes the proof
of lemma \ref{model lemma} and of theorem \ref{main theorem}.
$\Box$

\begin{remark} In the example of the decomposition $D^4=A'_2\cup B'_2$ in figures \ref{fig:A2},
\ref{fig:A2 Kirby} the proof above shows that $I_{\lambda}(A'_2,{\alpha})=1$. One may find an
explicit construction of a Bing cell bounded by ${\alpha}$ in $A'_2$ in the proof of
\cite[Lemma 7.3]{K1}.
\end{remark}

\section{Acknowledgements} This research was supported in part by the NSF under
grants DMS-0306934 and  NSF PHY05-51164.

\end{document}